\documentclass{imanum}
\usepackage{graphicx}
\usepackage{algorithm} 
\usepackage{algorithmic} 
\newcommand{\supp}{\operatorname{supp}}

\jno{drnxxx}

\begin{document}

\title{Wavelet Galerkin method for fractional elliptic differential
equations}

\author{%
{\sc
Weihua Deng\thanks{Corresponding author. Email: dengwh@lzu.edu.cn}, Yuwei Lin, Zhijiang Zhang} \\[2pt]
School of Mathematics and Statistics, Gansu Key Laboratory of Applied Mathematics and Complex Systems,
Lanzhou University, Lanzhou 730000, P.R. China
}

\maketitle

\begin{abstract}
{ Under the guidance of the general theory developed for classical partial differential equations (PDEs), we investigate the Riesz bases of wavelets in the spaces where fractional PDEs usually work, and their applications in numerically solving fractional elliptic differential equations (FEDEs). The technique issues are solved and the detailed algorithm descriptions are provided. Compared with the ordinary Galerkin methods, the wavelet Galerkin method we propose for FEDEs has the striking benefit of efficiency, since the condition numbers of the corresponding stiffness matrixes are small and uniformly bounded; and the Toeplitz structure of the matrix still can be used to reduce cost.  Numerical results and comparison with the ordinary Galerkin methods are presented to demonstrate the advantages of the wavelet Galerkin method we provide.
}
{fractional elliptic equation; Riesz bases; B-splines function; wavelet Garlerkin method; condition number
}

\end{abstract}

\section{Introduction}
\label{sec;introduction}
In recent decades, fractional operators have been playing more and more important roles in building the models [\cite{Diethelm:10}], e.g., in statistical physics (subdiffusion and superdiffusion), mechanics (theory of viscoelasticity and viscoplasticity), (bio-)chemistry (modelling of polymers and proteins), electrical engineering (transmission of ultrasound waves), medicine (modelling of human tissue under mechanical loads), etc. And in most of the cases, the models are appeared in the form of the fractional partial differential equations (PDEs), including the time dependent fractional PDEs and steady state fractional PDEs. Efficiently solving these fractional PDEs naturally becomes an urgent topic. Because of the nonlocal properties of fractional operators, obtaining the analytical solutions of the fractional PDEs is more challenging or sometimes even impossible; or the obtained analytical solutions are less valuable (expressed by transcendental functions or infinite series). Luckily, some important progress has been made for numerically solving the fractional PDEs by finite difference methods, e.g., see [\cite{Deng:14,Meerschaert:04,Sousa:12,Tian:14,Yuste:06,Zhuang:09}], finite element methods [\cite{Deng:08, Ervin:06}], spectral methods [\cite{Li:10, Zayernouri:13}], etc.

For the time dependent fractional PDEs, there are already some works to deal with the issue of computational efficiency, including the method of using the Toeplitz structure of the matrixes to reduce computational cost [\cite{Wang:12}] and the multigrid method [\cite{Chen:14, Pang:12}]. Heavy computational costs caused by the fast increasing of the condition numbers of the corresponding stiffness matrix with the mesh refinement and the inherent nonlocal properties of fractional operators are the main challenges that numerically solving the fractional elliptic differential equations (FEDEs) faces.  The condition number of the ordinary Galerkin equation is proportional to $h^{-2}$ for classical second order elliptical differential equations, and is proportional to $h^{-4}$ for fourth order ones [\cite{Jia:11}], where $h$ is the mesh size. And the condition number of the ordinary Galerkin equation for FEDEs with $\beta$-th order fractional derivative is proportional to $h^{-\beta}$  [\cite{Deng:13}]. So this problem is expected to be solved by the multiresolution methods including multigrid methods [\cite{Braess:95, Xu:90}] and wavelet methods [\cite{Roach:2000, Jia:11}]. It seems there are few works on the numerical methods for FEDEs [\cite{Ervin:06, Wang:13}], and almost no works for considering the efficiency of the numerical methods.

Wavelets have the strong multiresolution properties, and have been proven to be a powerful tool in signal and image processing such as image compression and denoising. In recent decades, the wavelet methods have also been well developed in solving the classical PDEs.  For the numerical treatment of PDEs, the efficiency of the wavelet method is greatly impacted by the properties of the wavelet bases; and in the sense of controlling the condition number, one can choose the Riesz bases of spline wavelets  [\cite{Jia:06, Jia:2006}]. Spline wavelets with short support are investigated in [\cite{Jia:2011}] and [\cite{Shen:06}]. The paper [\cite{Jia:09}] constructs the Riesz bases of spline wavelets on the interval $[0,1]$ with homogeneous boundary conditions. One can also refer to [\cite{Jia:11}] for the general theory of the construction of Riesz bases of wavelets and their applications to the numerical solutions of elliptic differential equations. Under the guidance of the theory being well developed for classical PDEs, in this paper we discuss the Riesz bases in the spaces where fractional PDEs usually work, and their applications in effectively solving FEDEs. The central gain of using the wavelet Galerkin method to solve FEDEs is its efficiency since the condition numbers of the corresponding stiffness matrix are small and uniformly bounded. The concrete FEDEs we discuss are the following one and two dimensional steady state fractional equations:
\begin{equation} \label{1.1}
-Da(p\,{} _0D_x^{-\beta} + q\,{} _x D_1^{-\beta})Du  = f(x),~~x \in \Omega=(0,1),
\end{equation}
with the boundary conditions $u(0)=u(1)=0$ and
\begin{equation}  \label{1.2}
-D_{x}^{s}a_{1}(p_{1}\,{} _0D_x^{-\alpha} + q_{1}\,{} _x D_1^{-\alpha})D_{x}u - D_{y}^{s}a_{2}(p_{2}\,{} _0D_y^{-\beta} + q_{2}\,{} _y D_1^{-\beta})D_{y}u = f(x,y),~~x,y \in \Omega=(0,1)^2,
\end{equation}
with $s=2$ or $s=3$; when $s=2$, the boundary conditions are   $u(x,y)|_{\partial \Omega}=0$, $(\partial u(x,y)/ \partial x)|_{x=0,y\in [0,1]}=0$, and $(\partial u(x,y)/ \partial y)|_{x\in [0,1],y=0}=0$; and when $s=3$, the boundary conditions are taken as $u(x,y)|_{\partial \Omega}=0$, $(\partial u(x,y)/ \partial x)|_{x=0~ {
\rm and}~ x=1,\,y\in [0,1]}=0$, and $(\partial u(x,y)/ \partial y)|_{x\in [0,1],\,y=0~ {\rm and} ~y=1}=0$. In (\ref{1.1}) and (\ref{1.2}), $a$, $a_1$, and $a_2$ are positive real numbers; $0\le\alpha,\,\beta<1$, $0 \le p,q,p_1,q_1,p_2,q_2 \le 1$ satisfying $p+q=p_1+q_1=p_2+q_2=1$; $D_x^s$ or $D_y^s$ means $s$ times partial derivative in $x$ or $y$ direction.
 The left and right Riemann-Liouville  fractional integral of the function $u(x)$ on $[a,b]$, $-\infty \leq a < b \leq \infty$, are respectively defined by [\cite{Podlubny:99}]
\begin{equation}\label{1.3}
 _{a}D_x^{-\alpha}u(x)=
\frac{1}{\Gamma(\alpha)} \int_{a}\nolimits^x{\left(x-\xi\right)^{\alpha-1}}{u(\xi)}d\xi,
\end{equation}
and
\begin{equation}\label{1.4}
 _{x}D_{b}^{-\alpha}u(x)=
 \frac{1}{\Gamma(\alpha)} \int_{x}\nolimits^{b}{\left(\xi-x\right)^{\alpha-1}}{u(\xi)}d\xi.
\end{equation}

The  outline of this paper is as follows. In Section 2, the fractional Sobolev space and a class of B-spline functions are firstly introduced; based on these functions, we introduce the Riesz bases in one and two dimensional fractional Sobolev spaces. In Section 3, we discuss the wavelet Galerkin method for FEDEs, present its detailed algorithm description, and the extensive numerical experiments are also performed to show its powerfulness. We conclude the paper with some remarks in the last section.

\section{Wavelet Riesz bases in fractional Sobolev space}\label{sec:1}
The Riesz bases play a vital role in controlling the condition number of stiffness matrix when using wavelet Galerkin method to solve the FEDEs. We present the Riesz bases in one and two dimensional fractional Sobolev spaces. First, we introduce the spaces where FEDEs work and their relations to the fractional order Hilbert spaces $H_{0}^{\mu}$.


\subsection{Fractional Sobolev space}
We introduce the abstract setting for FEDEs, including the left, right, and symmetric fractional derivative spaces; and then show the equivalence of the fractional derivative spaces with fractional order Hilbert spaces [\cite{Ervin:06}].


\begin{definition}[left fractional derivative] \label{DLFD} Let $u$ be a function defined on
$R$, $\mu>0$, $n$ be the smallest integer greater than $\mu$
$(n-1 \le \mu <n)$, and $\sigma=n-\mu$. Then the left fractional
derivative of order $\mu$ is defined to be
$$
{\bf D}^\mu u:=D^n {_{-\infty}D_x^{-\sigma}}
u(x)=\frac{1}{\Gamma(\sigma)} \frac{d^n}{dx^n}
\int_{-\infty}^x(x-\xi)^{\sigma-1}u(\xi)d\xi.
$$
\end{definition}

\begin{definition}[right fractional derivative] \label{DRFD} Let $u$ be a function defined on
$R$, $\mu>0$, $n$ be the smallest integer greater than $\mu$
$(n-1 \le \mu <n)$, and $\sigma=n-\mu$. Then the right
fractional derivative of order $\mu$ is defined to be
$$
{\bf D}^{\mu\ast} u:=(-D)^n {_xD_{\infty}^{-\sigma}}
u(x)=\frac{(-1)^n}{\Gamma(\sigma)} \frac{d^n}{dx^n}
\int_x^\infty(\xi-x)^{\sigma-1}u(\xi)d\xi.
$$
\end{definition}
{\em Note:} If $\supp(u) \subset (a,b)$, then ${\bf D}^{\mu}
u={_aD_x^\mu} u$ and ${\bf D}^{\mu\ast} u={_xD_b^\mu} u$,
where ${_aD_x^\mu} u$ and ${_xD_b^\mu} u$ are the left and right
Riemann-Liouville fractional derivative of order $\mu$ defined as
$$
{_aD_x^\mu} u=\frac{1}{\Gamma(\sigma)} \frac{d^n}{dx^n}
\int_a^x(x-\xi)^{\sigma-1}u(\xi)d\xi,
$$
and
$$
{_xD_b^\mu} u=\frac{(-1)^n}{\Gamma(\sigma)} \frac{d^n}{dx^n}
\int_x^b(\xi-x)^{\sigma-1}u(\xi)d\xi.
$$

\begin{definition}[left fractional derivative space]
Let $\mu>0$. Define the semi-norm
\begin{displaymath}
\mid u \mid_{J_{L}^{\mu}(R)} := \parallel {\bf D}^{\mu}u \parallel_{L^{2}(R),}
\end{displaymath}
and norm
\begin{displaymath}
\parallel u \parallel_{J_{L}^{\mu}(R)} := (\parallel u \parallel_{L^{2}(R)}^{2} + \mid u \mid_{J_{L}^{\mu}(R)}^{2})^{1/2},
\end{displaymath}
and let $J_{L}^{\mu}(R)$ denote the closure of $C_{0}^{\infty}(R)$ with respect to $\parallel \cdot \parallel_{J_{L}^{\mu}(R)}$.
\end{definition}

\begin{definition}[right fractional derivative space]
Let $\mu>0$. Define the semi-norm
\begin{displaymath}
\mid u \mid_{J_{R}^{\mu}(R)} := \parallel {\bf D}^{\mu*}u \parallel_{L^{2}(R),}
\end{displaymath}
and norm
\begin{displaymath}
\parallel u \parallel_{J_{R}^{\mu}(R)} := (\parallel u \parallel_{L^{2}(R)}^{2} + \mid u \mid_{J_{R}^{\mu}(R)}^{2})^{1/2},
\end{displaymath}
and let $J_{R}^{\mu}(R)$ denote the closure of $C_{0}^{\infty}(R)$ with respect to $\parallel \cdot \parallel_{J_{R}^{\mu}(R)}$.
\end{definition}

The Fourier transform of a function $f \in L^1(R)$ is defined by
\begin{displaymath}
\hat{f}(\xi)：=\frac{1}{\sqrt{2\pi}}\int\nolimits_{R} f(x)e^{-ix \cdot \xi} \,\mathrm{d}x , ~~~\xi \in R.
\end{displaymath}
The Fourier transform can be naturally extended to functions in $L^2(R)$. For $\mu > 0$, we denote by $H^{\mu}(R)$ the Sobolev space of all functions $f \in L^2(R)$ such that the seminorm
\begin{equation}
|f|_{H^{\mu}(R)}:=\left(\frac{1}{2\pi}\int\nolimits_{R}|\hat{f}(\xi)|^{2}|\xi|^{2\mu}\,\mathrm{d}\xi\right)^{1/2}
\end{equation}
is finite. The space $H^{\mu}(R)$ is a Hilbert space with the inner product given by
\begin{displaymath}
\langle f,g\rangle_{H^{\mu}(R)}:=\frac{1}{2\pi}\int\nolimits_{R}\hat{f}(\xi)\overline{\hat{g}(\xi)}[1+|\xi|^{2\mu}]\,\mathrm{d}\xi,~~~f,g \in H^{\mu}(R).
\end{displaymath}
The corresponding norm in $H^{\mu}(R)$ is given by $\parallel f \parallel_{H^{\mu}(R)} := \sqrt{\parallel f \parallel^{2}_{L^{2}(R)}+|f|^{2}_{H^{\mu}(R)}}$.

\begin{lemma}[\cite{Ervin:06}]
Let $\mu>0$. The spaces $J_{L}^{\mu}(R)$, $J_{R}^{\mu}(R)$, and $H^{\mu}(R)$ are equal with equivalent semi-norms and norms.
\end{lemma}

%
%

\begin{definition}[symmetric fractional derivative space]
Let $\mu >0$, $\mu \neq n - 1/2$, $n \in N$. Define the semi-norm
\begin{displaymath}
\mid u \mid_{J_{S}^{\mu}(R)} := \mid ({\bf D}^{\mu}u, {\bf D}^{\mu*}u) \mid_{L^{2}(R)}^{1/2},
\end{displaymath}
and norm
\begin{displaymath}
\parallel u \parallel_{J_{S}^{\mu}(R)} := (\parallel u \parallel_{L^{2}(R)}^{2} + \mid u \mid_{J_{S}^{\mu}(R)}^{2})^{1/2},
\end{displaymath}
and let $J_{S}^{\mu}(R)$ denote the closure of $C_{0}^{\infty}(R)$ with respect to $\parallel \cdot \parallel_{J_{S}^{\mu}(R)}$.
\end{definition}

\begin{lemma}[\cite{Ervin:06}]
For $\mu>0$, $\mu \neq n - 1/2$, $n \in N$, the spaces $J_{L}^{\mu}(R)$ and $J_{S}^{\mu}(R)$ are equal, with equivalent semi-norms and norms.
\end{lemma}

\begin{definition}
Define the spaces $J_{L,0}^{\mu}(\Omega)$, $J_{R,0}^{\mu}(\Omega)$, $J_{S,0}^{\mu}(\Omega)$, and $H_{0}^{\mu}(\Omega)$ as the closures of $C_{0}^{\infty}(\Omega)$ under their respective norms.
\end{definition}

We next turn to the equivalence of the fractional derivative spaces $J_{L,0}^{\mu}(\Omega)$, $J_{R,0}^{\mu}(\Omega)$, $J_{S,0}^{\mu}(\Omega)$, and the fractional order Hilbert space $H_{0}^{\mu}(\Omega)$.

%

\begin{lemma}[\cite{Ervin:06}]
 Let $\mu>0$. Then the spaces $J_{L,0}^{\mu}(\Omega)$, $J_{R,0}^{\mu}(\Omega)$, and $H_{0}^{\mu}(\Omega)$ are equal. Also, if $\mu \neq n - 1/2$, $n \in N$, the spaces $J_{L,0}^{\mu}(\Omega)$, $J_{R,0}^{\mu}(\Omega)$, and $H_{0}^{\mu}(\Omega)$ have equivalent semi-norms and norms.
\end{lemma}

\subsection{Wavelet bases and the related lemmas}
 Let $N$ denote the set of positive integers, $J$ be a (finite or infinite) countable set. By $\ell(J)$ we denote the linear space of all complex-valued sequences $(u_{j})_{j \in J}$; $\ell_{0}(J)$ denotes the linear space of all sequences $(u_{j})_{j \in J}$ with only finite nonzero terms; and $\ell^{2}(J)$ denotes the linear space of all sequences $u=(u_{j})_{j \in J}$ such that $\parallel u \parallel_{2} := (\sum_{i \in J}|u_{j}|^{2})^{1/2} < \infty$.



Let $H$ be a Hilbert space. A sequence $\{v_j\}_{j\in J}$ in $H$ is said to be a \textbf{ Riesz sequence} if there exist two positive constants $C_{1}$ and $C_{2}$ such that the inequalities
\begin{equation}\label{equ:1.1.1}
C_{1}(\sum_{j \in J}|c_{j}|^{2})^{1/2} \leq \| \sum_{j \in J} c_{j}v_{j} \| \leq C_{2}(\sum_{j \in J}|c_{j}|^{2})^{1/2}
\end{equation}
hold true for every sequence $(c_{j})_{j \in J}$ in $\ell_{0}(J)$. If this is the case, then the series $\sum_{j \in J}c_{j}v_{j}$ converges unconditionally for every $(c_{j})_{j \in J}$ in $\ell^{2}(J)$, and the inequalities in ~\eqref{equ:1.1.1} are valid for all $(c_{j})_{j \in J}$ in $\ell^{2}(J)$. We call $C_{1}$ a Riesz lower bound and $C_{2}$ a Riesz upper bound. If $\{v_j\}_{j\in J}$ is a Riesz sequence in $H$, and the linear span of $\{v_j\}_{j\in J}$  is dense in $H$, then $\{v_j\}_{j\in J}$  is a \textbf{Riesz basis} of $H$.

In numerical simulations, spline wavelet bases are more popular, since they are relatively smooth, have a small support, and can be got in a closed form. The widely and effectively way to build the Riesz bases is based on multiresolution analysis. Following [\cite{Jia:09}], we first introduce the Riesz bases in one dimension, then extend them to two dimensional case.

\subsubsection{Riesz ~bases ~in ~$H_{0}^{\mu}((0,1))$ }
For a positive integer $m$, let $M_{m}$ denote the $B$-spline of order $m$, which is the convolution of $m$ copies of the characteristic function of the interval $[0,1]$: 
\begin{displaymath}
M_{m}(x) = \int_{0}^{1} M_{m-1}(x-t)\,\mathrm{d}t,~~~x \in R,
\end{displaymath}
where $M_{1} :=\chi_{[0,1]}$. For $m=2$ and $3$, the spline functions are given as follows,
\[
  M_{2}(x) = \left\{
        \begin{array}{ll}
            x,  &\mbox{$0 \leq x \leq 1$,}\\
            2-x,  &\mbox{$1 \leq x \leq 2$.}\\
           \; 0,    & \quad else;
        \end{array}
        \right.\quad\quad
\]
\[
  M_{3}(x) = \left\{
        \begin{array}{lll}
            \frac{1}{2}x^{2},  &\mbox{$0 \leq x \leq 1$,}\\
            -x^{2} + 3x - \frac{3}{2},  &\mbox{$1 \leq x \leq 2$,}\\
            \frac{1}{2}x^{2} - 3x + \frac{9}{2},  &\mbox{$2 \leq x \leq 3$.}\\
            \;0,                   & \quad else,

        \end{array}
        \right.
\]
and they, respectively, satisfy the following refinement equations,
\begin{equation}
M_{2}(x) = \frac{1}{2}M_{2}(2x) + M_{2}(2x-1) + \frac{1}{2}M_{2}(2x-2);
\end{equation}
and
\begin{equation}
M_{3}(x) = \frac{1}{4}M_{3}(2x) + \frac{3}{4}M_{3}(2x-1) + \frac{3}{4}M_{3}(2x-2) +\frac{1}{4}M_{3}(2x-3).
\end{equation}


From the definition, it follows immediately that $M_{m}$ is supported on $[0,m]$, $M_{m}(x) > 0$ and $M_{m}(m-x) = M_{m}(x)$ for $0 < x < m$. Moreover, $M_{m} \in H_{0}^{\mu}(0,m)$ for $0 < \mu < m - 1/2$. Let
\begin{eqnarray*}
\phi_{n,j}(x) &:=& 2^{n/2}M_{r}(2^{n}x - j), \quad j \in I_{n} := \{0,1,\ldots,2^{n}-r\}.
\end{eqnarray*}
 Then there exists $n_0\in N$, such that $n \geq n_{0},\,V_{n} := {\rm span}\{\phi_{n,j}:j \in I_{n}\}$ is a subspace of $H_{0}^{\mu}(0,1)$ for $0 \leq \mu \leq r-1/2$. Evidently, $V_{n} \subset V_{n+1}$, for $r \geq 2$, each function $f$ in $V_{n}$ satisfies the homogeneous boundary conditions
\begin{displaymath}
f^{(k)}(0) = f^{(k)}(1) = 0,   ~~~k = 0,1,\ldots,r-2.
\end{displaymath}
And the following result on approximation accuracy holds:
\begin{equation}
\inf_{v\in V_n}\|u-v\|_{H^\mu(0,1)} \leq C(1/2^{n})^{r-\mu}|u|_{H^{r}(0,1)},    ~~~u \in  H_0^{r-1}(0,1)\cap H^{r}(0,1).
\end{equation}
Its proof is similar to the one given in [\cite{Jia:11}].

For the construction of wavelet bases, suppose that $t\in N$, $r \geq t$, and $r+t$ is an even integer. Let $n_{0}$ be the least integer such that $2^{n_{0}} \geq r + t$, and define $\hat{\phi}_{n,j}(x)$
\begin{eqnarray*}
  \hat{\phi}_{n,j}(x) &:=& 2^{n/2}M_t(2^{n}x-j-(r-t)/2).
\end{eqnarray*}
Let $\hat{V}_{n} := {\rm span} \{\hat{\phi}_{n,j}:j \in I_{n}\}$, obviously, $\hat{\phi}_{n,j}(x) = 0$ for $x \in R\setminus[0,1]$ when $n \geq n_{0}$, $\hat{V}_{n} \subset \hat{V}_{n+1}$. Then we find the direct sum decomposition of $ V_{n+1}$ ($V_n\bigoplus W_n$) and $\hat{V}_{n+1}$ ($\hat{V}_n\bigoplus \hat{W}_n$) by demanding that $W_{n} := V_{n+1} \cap \hat{V}_{n}^{\bot}$ and $\hat{W}_{n} := \hat{V}_{n+1} \cap V_{n}^{\bot}$, respectively. The desired wavelet bases for $W_n$ and $\hat{W}_{n}$ can be constructed by studying the slant matrixes.
Two important wavelet bases for $r=2$ and $3$ are given as follows:

~$(a)$ For $r=2$ and $t=2$, let
\begin{equation}
\psi(x) = \frac{1}{24}M_{2}(2x) - \frac{1}{4}M_{2}(2x-1) + \frac{5}{12}M_{2}(2x-2) - \frac{1}{4}M_{2}(2x-3) + \frac{1}{24}M_{2}(2x-4),
\end{equation}
and
\begin{equation}
\psi_{1}(x) = \frac{3}{8}M_{2}(2x) - \frac{1}{4}M_{2}(2x-1) + \frac{1}{24}M_{2}(2x-2).
\end{equation}
For $n \geq 2$ and $x \in R$, we define
\[
  \psi_{n,j}(x):=\left\{
            \begin{array}{lll}
               2^{n/2}\psi_{j}(2^{n}x),  &\mbox{$j=1$,}\\
                2^{n/2}\psi(2^{n}x-j+2),  &\mbox{$j=2,\ldots,2^{n}-1$,}\\
                2^{n/2}\psi_{2^{n}-j+1}(2^{n}(1-x)),  &\mbox{$j=2^{n}$.}
            \end{array}
            \right.
\]

~$(b)$ For $r=3$ and $t=1$, let
\begin{equation}
\psi(x) = \frac{1}{12}M_{3}(2x) - \frac{5}{12}M_{3}(2x-1) + \frac{5}{12}M_{3}(2x-2) - \frac{1}{12}M_{3}(2x-3),
\end{equation}
and
\begin{equation}
\psi_{1}(x) = \frac{5}{12}M_{3}(2x) - \frac{1}{12}M_{3}(2x-1).
\end{equation}
For $n \geq 2$ and $x \in R$, we define
\[
  \psi_{n,j}(x):=\left\{
            \begin{array}{lll}
               2^{n/2}\psi_{j}(2^{n}x),  &\mbox{$j=1$,}\\
                2^{n/2}\psi(2^{n}x-j+2),  &\mbox{$j=2,\ldots,2^{n}-1$,}\\
                2^{n/2}\psi_{2^{n}-j+1}(2^{n}(1-x)),  &\mbox{$j=2^{n}$.}
            \end{array}
            \right.
\]
Then we have the following important lemma.
\begin{lemma}[\cite{Jia:09}]\label{lemma:3.1}
For $n \geq n_{0}$ and $j \in J_{n} := \{1,2,\ldots,2^{n}\}$, let $\psi_{n,j}$ be the functions as constructed above. Then the set
\begin{displaymath}
\{2^{-n_{0}\mu}\phi_{n_{0},j} : j \in I_{n_{0}}\} \cup \bigcup_{n=n_{0}}^{\infty}\{2^{-n\mu}\psi_{n,j} : n \geq n_{0}, j \in J_{n}\},
\end{displaymath}
forms a Riesz basis of $H_{0}^{\mu}(0,1)$ for $0 < \mu < r-\frac{1}{2}$.
\end{lemma}

\subsubsection{Riesz ~bases ~in ~$H_{0}^{\mu}((0,1)^2)$ }

In order to obtain a Riesz basis of ~$H_{0}^{\mu}((0,1)^2)$, we use the tensor product denoted by $\otimes$. For two functions $v$ and $\omega$ defined on $(0,1)$, we use $v \otimes \omega$ to denote the function on $(0,1)^{2}$ given by
\begin{displaymath}
v\otimes\omega(x,y) := v(x)\omega(y), ~~~0 \leq x,y \leq 1.
\end{displaymath}
For $n \geq n_0, n\in N$, let $I_n:=\{j=(j_{1},j_{2})\in Z^{2}:0 \leq j_{1} \leq 2^{n}-r,\,0 \leq j_{2} \leq 2^{n}-r\}$. We denote the approximate space of $H_{0}^{\mu}((0,1)^{2})$ by $(\tilde{V}_{n})_{n \geq n_{0}}$. Define
\begin{eqnarray}
&& \widetilde{\phi}_{n,j} ~:=~ \phi_{n,j_{1}} \otimes \phi_{n,j_{2}}, ~j \in I_n;\nonumber \\
&& \widetilde{\Phi}_{n} ~:=~ \{\widetilde{\phi}_{(n,j},j=(j_{1},j_{2}) \in I_{n}\};\\
&& \widetilde{V}_{n} ~:=~ span\{\widetilde{\Phi}_{n}\}.\nonumber
\end{eqnarray}
Similarly we can define the corresponding $\widetilde{\hat{\phi}}_{n,j}$.

For the sequence of the subspaces $\widetilde{V}_{n}$, we have the following properties
\begin{eqnarray*}
&\bullet& \widetilde{V}_{n_0} \subset \widetilde{V}_{n_0+1} \subset \widetilde{V}_{n_0+2} \subset \ldots ;\\
&\bullet& \bigcup_{n=n_0}^{\infty}\widetilde{V}_{n}~is~dense~in~H_{0}^{\mu}((0,1)^2)~for~0 < \mu < r-1/2;\\
&\bullet& dim(\widetilde{V}_{n})=(2^{n}-r+1)^2;
\end{eqnarray*}



Furthermore, define
\begin{eqnarray*}
\Gamma_{n}^{\prime} &:=& \{\phi_{n,j_{1}} \otimes \psi_{n,j_{2}}:0 \leq j_{1} \leq 2^{n}-r, 1 \leq j_{2} \leq 2^{n}\};\\
\Gamma_{n}^{\prime\prime} &:=& \{\psi_{n,j_{1}} \otimes \phi_{n,j_{2}}: 0 \leq j_{2} \leq 2^{n}-r,1 \leq j_{1} \leq 2^{n}\};\\
\Gamma_{n}^{\prime\prime\prime} &:=& \{\psi_{n,j_{1}} \otimes \psi_{n,j_{2}}:1 \leq j_{1} \leq 2^{n},1 \leq j_{2} \leq 2^{n}\}.
\end{eqnarray*}


For $n\geq n_{0}$, let $\Gamma_{n} := \Gamma_{n}^{\prime} \bigcup \Gamma_{n}^{\prime\prime} \bigcup \Gamma_{n}^{\prime\prime\prime}$ ,
and $\widetilde{W_{n}}~:=~ span\{\Gamma_{n} \}$, then $\Gamma_{n}$ is a Riesz basis of ~$\widetilde{W_{n}}$ in the $L_{2}$ space.
The dimensions satisfied the following relation:
\begin{eqnarray*}
dim(\tilde{V}_{n+1}) &=& dim(\tilde{V}_{n})+dim(\tilde{W}_{n})\\
&=& (2^{n}-r+1)^2+(2^{n}-r+1)2^{n}+2^{n}(2^{n}-r+1)+2^{2n}\\
&=& (2^{n+1}-r+1)^2.
\end{eqnarray*}
For every $f\in L^2(0,1)^2$ and $n\geq n_0$, let $P_n f$ be the unique element in $\widetilde{V}_{n}$ such that
\begin{displaymath}
<P_nf,\widetilde{\hat{\phi}_{n,j}}>=<f,\widetilde{\hat{\phi}_{n,j}}>, \quad \forall j\in I_n.
\end{displaymath}
It is easy to check that $P_n$ is a projector from $L^2(0,1)^2$ onto $\widetilde{V}_{n}$, $\widetilde{W_{k}}$ is the kernel space of $P_n$, and $\widetilde{V}_{n+1}$ is the direct sum of $\widetilde{V}_{n}$ and $\widetilde{W_{n}}$. Using the similar way of the proof to the one dimensional case given in [\cite{Jia:09}], we can prove the following lemma.
\begin{lemma}\label{lemma:3.2}
For ~$0 < \mu < r-1/2$, the set
\begin{displaymath}
\{2^{-n_{0}\mu}\widetilde{\Phi}_{n_{0}}\} \cup \bigcup_{n=n_{0}}^{\infty}\{2^{-n\mu}\Gamma_{n})
\end{displaymath}
forms a Riesz basis of ~$H_{0}^{\mu}((0,1)^{2})$.
\end{lemma}

We have discussed the Riesz bases for the fractional Sobolev spaces in one and two dimensional cases. Lemma \ref{lemma:3.3} will present the reason of constructing the Riesz bases and give the theoretical foundation of the following numerical analysis.

\begin{lemma}[\cite{Jia:11}]\label{lemma:3.3}
If $\Psi^{norm}$ is a Riesz basis of $H_{0}^{\mu}(\Omega)$ and the bilinear form $a(u,v)$ corresponding to the equation is symmetric, continuous, and coercive, then the condition number of the stiffness matrix associated with $\Psi_{n}^{norm}$, i.e., $B_{n}:=(a(\chi,\psi))_{\chi,\psi \in \Psi_{n}^{norm}}$, is uniformly bounded.
\end{lemma}

\section{Wavelet Galerkin method for fractional elliptic differential equation}\label{sec:1}
In this section, we apply the Riesz bases to FEDEs in one and two dimensional spaces and present the corresponding algorithms. The provided methods are compared with the ordinary Galerkin method from the condition number and computational time that show the efficiency of wavelet Galerkin method.

\subsection{Wavelet Galerkin method for one dimensional FEDE}
For the one dimensional space, we have provided the Riesz bases of the fractional Sobolev spaces in Lemma \ref{lemma:3.1}. And they can be used to solve the following eqution
\begin{equation}\label{3.1}
-Da(p\, {}_0D_x^{-\beta} + q\,{} _x D_1^{-\beta})Du  = f,~~x \in \Omega=(0,1),
\end{equation}
where $D$ represents the first spatial derivative, ${} _0D_x^{-\beta}$ and ~${} _xD_1^{-\beta}$ are the left and right Riemann-Liouville fractional integral operators, respectively, with $0 \leq \beta < 1$ and $p + q = 1$.


The variational formulation of this fractional differential equation is as follows
\begin{equation}
B(u,v)=(f,v),  ~~ \forall v \in H_0^\mu(0,1),
\end{equation}
where $B(u,v)=ap\langle {} _0D_x^{-\beta}Du,Dv \rangle + aq\langle {} _x D_1^{-\beta}Du,Dv \rangle$, $\mu = \frac{2-\beta}{2}$, and $\frac{1}{2} < \mu \leq 1$. According to [\cite{Ervin:06}], the above variational formulation has the unique solution in space $H_0^{\mu}(0,1)$. In order to numerically solve the resulting  variational formulation, we can use the subspace $V_{n}$ to approximate the space $H_{0}^{\mu}(0,1)$, i.e., find a $u_{n} \in V_{n}$ such that
\begin{equation}\label{equ:3.1.3}
ap\langle {} _0D_x^{-\beta}Du_{n},Dv \rangle + aq\langle {} _x D_1^{-\beta}Du_{n},Dv \rangle = \langle f,v \rangle ,  ~~\forall v \in V_{n}.
\end{equation}
Suppose that $\Phi_{n} := \{\phi_{n,j}:j \in I_n\}$ being a basis of $V_{n}$ and $u_{n} = \sum\nolimits_{\phi \in \Phi_{n}}y_{\phi}\phi$. Let $A_{n}$ be the matrix $(ap\langle {} _0D_x^{-\beta}D\sigma,D\phi \rangle +  aq\langle {} _xD_1^{-\beta}D\sigma,D\phi \rangle)_{\sigma,\phi \in \Phi_{n}}$, and $\xi_{n}$ the column vector $(\langle f,\phi \rangle)_{\phi \in \Phi_{n}}$. Then the column vector $y_{n}=(y_{\phi})_{\phi \in \Phi_{n}}$ is the solution of the system of linear equations
\begin{equation}\label{equ:3.1.4}
A_{n}y_{n} = \xi_{n}.
\end{equation}

We also employ the Riesz bases constructed in the previous section to solve the variational problem. For $n \geq n_{0}$, we have $\Psi_{n} := \{2^{-n_{0}\mu}\phi_{n_{0},j}:j \in I_{n_{0}}\} \cup \bigcup_{k=n_{0}}^{n-1}\{2^{-k\mu}\psi_{k,j}:j \in J_{k}\}$. Similarly, find a column vector $z_{n}=(z_{\psi})_{\psi \in \Psi_{n}}$ to satisfy
\begin{equation}\label{equ:3.1.5}
B_{n}z_{n} = \eta_{n},
\end{equation}
where $B_{n}=(ap\langle {} _0D_x^{-\beta}D\chi,D\psi \rangle + aq\langle {} _xD_1^{-\beta}D\chi,D\psi \rangle)_{\chi,\psi \in \Psi_{n}}$; and $\eta_{n}$ denote the column vector $(\langle f,\psi \rangle)_{\psi \in \Psi_{n}}$. Hence, $u_{n}=\sum_{\psi \in \Psi_{n}}z_{\psi}\psi$ is the approximate solution of $u$ in $V_{n}$.

Since $\Phi_{n}$ and $\Psi_{n}$ are two different bases of $V_{n}$, there is a unique transformation $S_{n}$ between $\Phi_{n}$ and $\Psi_{n}$, such that $\Psi_{n}$=$S_{n}\Phi_{n}$, which is called wavelet transformation. So we have
$
B_{n} = S_{n}A_{n}S_{n}^{T}
$
,
$
\eta_{n} = S_{n}\xi_{n}
$
and that \eqref{equ:3.1.4} is equivalent to \eqref{equ:3.1.5}. If we set $y_{n} = S_{n}^{T}z_{n}$, then \eqref{equ:3.1.4} is preconditioned by the transformation $S_{n}$. And $S_{n}$ is called a preconditioner.

Now, we present the algorithm of generating the matrixes $A_n$ and $B_{n}$ in Algorithm 1 (for performing the numerical analysis, sometimes we need explicitly to get the matrixes). For the corresponding classical problems, the matrix $A_{n}$ is symmetric and sparse. But for the fractional problems, $A_{n}$ is dense and also nonsymmetric when $p \neq q$. Fortunately, because $\phi_{n,j}$ are the dilation and translation of one single function $M_r(x)$, the matrix $A_{n}$ has the Toeplitz (diagonal-constant) structure. Then we only need to produce the first row and column of $A_n$, which greatly reduces the computation and storage costs.
\begin{algorithm}[!h t b]\selectfont \label{Alg:1}
\caption{Generating matrix $A_{n}$ and $B_n$ for $1D$}
\begin{algorithmic}[1]
\FOR{$i=0,1,\ldots, r-1$}
\STATE $a_1(i)=\langle {}_0D_x^{-\beta}\phi_{n,i}',\phi_{n,0}'\rangle$
\ENDFOR
\FOR{$j=0,1,\ldots,2^n-r$}
\STATE $q_1(j)=\langle {}_0D_x^{-\beta}\phi_{n,0}',\phi_{n,j}'\rangle$
\ENDFOR
\STATE Initialize the unit matrix $P$: $P=speye(2^{n_0}-r+1)$
\FOR{$k = n_0+1,\ldots,n$}
\STATE Let $P_{k}$ satisfy: \\[4pt]
 \hspace{16pt} ${\Phi_{k-1} \choose 2^{-(k-1)\mu}\Gamma_{k-1}} = P_k\Phi_{k}$
\STATE $ P_{k}(1:2^{k-1}-r+1,1:2^{k-1}-r+1) := P \times P_{k}(1:2^{k-1}-r+1,1:2^{k-1}-r+1)$\\
$P_{k}(1:2^{k-1}-r+1,2^{k-1}-r+2:end) := P \times P_k(1:2^{k-1}-r+1,2^{k-1}-r+2:end)$\\
$ P:=P_k$
\ENDFOR
\STATE $P(1:2^{n_0}-r+1,:)=2^{-n_0\mu}P(1:2^{n_0}-r+1,:)$
\STATE Get the transformation matrix $S_n$: $S_n=P$
\STATE $A_n=a*p\mathcal{T}(q_1,a_1)+a*q\mathcal{T}(a_1,q_1)$
\STATE $B_n:=S_nA_nS_n^T$\\[8pt]
\end{algorithmic}
\begin{algorithmic}
\STATE{\em Note:} $\Gamma_k:=\{\psi_{k,j}:j \in J_{k}\}$, and $\mathcal{T}(col,row)$ denote the Toeplitz matrix produced by its first column $col$ and its first row $row$, $\mu=1-\beta/2$.
\end{algorithmic}
\end{algorithm}

When taking $p=q=0.5$ in (\ref{3.1}), the matrixes $A_{n}$ and $B_{n}$ are both symmetrical; in Table 1, it is shown that the increasing of the condition numbers of $A_{n}$ is as $\mathcal{O}(1/h^{2-\beta})$, where $h$ is the mesh size; and the condition numbers of the corresponding $B_{n}$ are uniformly bounded w.r.t $h$, which confirms Lemma \ref{lemma:3.3}. The observations also hold for the nonsymmetrical case with $p=1$ and $q=0$ in (\ref{3.1}), see Table 2.
\begin{table}[!h t b]\fontsize{7.0pt}{12pt}\selectfont \label{Tab1}
  \begin{center}
  \caption {The condition numbers of the matrixes $A_{n}$ and $B_{n}$ with $p = 0.5$, $q = 0.5$, and $r = 2$} \vspace{5pt}
  \begin{tabular*}{\linewidth}{@{\extracolsep{\fill}}*{6}{c}}                                    \hline  
$\beta$ &$n$ & size of $B_{n}$  & $\kappa(A_{n})$        & $\kappa(B_{n})$ &      \\\hline
              &         3&    7 $\times$ 7&     10.0502&     8.7751  \\
              &         4&    15 $\times$ 15&     28.4815&     10.0365 \\
              &         5&    31 $\times$ 31&     80.6947&     10.6426  \\
$\beta = 0.5$ &         6&     63 $\times$ 63&     228.5288&     11.0895   \\
              &         7&     127 $\times$ 127&     646.8779&     11.4778   \\
              &         8&     255 $\times$ 255&     1.8304e+03&     11.8235 \\
              &         9&     511 $\times$ 511&     5.1784e+03&     12.1302  \\
              &       10&     1023$\times$1023&       1.4648e+04&     12.4007   \\\hline
               &        3&    7 $\times$ 7&     6.2382&     9.2001  \\
               &        4&    15 $\times$ 15&     14.7486&     10.1688  \\
               &        5&    31 $\times$ 31&     35.0981&     10.7312 \\
$\beta = 0.75$ &        6&     63 $\times$ 63&     83.6018&     11.2143   \\
               &        7&     127 $\times$ 127&     199.0556&     11.6654   \\
               &        8&     255 $\times$ 255&     473.7381&     12.0813 \\
               &        9&     511 $\times$ 511&     1.1271e+03&    12.4573 \\
               &       10&     1023$\times$1023&     2.6813e+03 &      12.7929 \\\hline
    \end{tabular*}
  \end{center}
\end{table}

\begin{table}[!h t b]\fontsize{7.0pt}{12pt}\selectfont
  \begin{center}
  \caption {The condition numbers of the matrixes $A_{n}$ and $B_{n}$ with $p = 1$, $q = 0$, $r = 2$} \vspace{5pt}
  \begin{tabular*}{\linewidth}{@{\extracolsep{\fill}}*{6}{c}}                                    \hline  
$\beta$ &$n$ & size of $B_{n}$  & $\kappa(A_{n})$        & $\kappa(B_{n})$      \\\hline
              &         3&    7 $\times$ 7&     8.3362&     6.6338  \\
              &         4&    15 $\times$ 15&     23.2013&     7.6409  \\
              &         5&    31 $\times$ 31&     65.3566&     8.7345  \\
$\beta = 0.5$ &         6&     63 $\times$ 63&     184.6258&     9.5570   \\
              &         7&     127 $\times$ 127&     522.0054&     10.2252   \\
              &         8&     255 $\times$ 255&     1.4763e+03&     10.7896 \\
              &         9&     511 $\times$ 511&     4.1754e+03&     11.2744  \\
              &        10&     1023$\times$1023&      1.1810e+04&      11.6965  \\\hline

               &        3&    7 $\times$ 7&     6.2421&     6.8589  \\
               &        4&    15 $\times$ 15&     14.2077&     8.5584 \\
               &        5&    31 $\times$ 31&     33.2186&     9.8513  \\
$\beta = 0.75$ &        6&     63 $\times$ 63&     78.4138&     10.6949   \\
               &        7&     127 $\times$ 127&     185.1202&     11.9194   \\
               &        8&     255 $\times$ 255&     441.7010&      12.7937 \\
               &        9&     511 $\times$ 511&     1.0501e+03&     13.5884  \\
               &       10&     1023$\times$1023&       2.4971e+03&      14.3144\\\hline

\end{tabular*}
\end{center}
\end{table}
For further showing the powerfulness of the wavelet Galerkin method (solving the algebraic equation w.r.t. $B_n$), we use the Krylov subspace method to solve the algebraic equations w.r.t. $A_n$ and $B_n$, respectively. In fact, solving the algebraic system of $B_n$ is essentially to solve the preconditioned system of $A_n$. It is well known that the conjugate gradient method has the properties of short recursive and residuals minimality; a small condition number usually means a fast iterative speed; and it can only  be used to the symmetric positive define systems. For the nonsymmetrical
system, the Bi-CGSTAB method is popular, since it retains the property of short recurrence, usually have a fast convergence speed compared with the other Krylov subspace methods; but the interruption may occurs.  The algorithm of wavelet preconditioned Bi-CGSTAB method is given in Algorithm 2. For computing $A_nx$, we use the method of Toeplitz fast matrix-vector multiplications proposed in Algorithm 3, with the computational complexity just as $\mathcal{O}(n2^n)$. And for $S_nx$ and $S_n^Tx$, the fast wavelet transformation (FWT) or the sparsity of $S_n$  can be applied, which just has the computational complexity as $\mathcal{O}(2^n)$.
\begin{algorithm}[h t b p]\selectfont \label{Alg:2}
\caption{Wavelet preconditioned Bi-CGSTAB}
\begin{algorithmic}[1]
\STATE Given the initial value of interation $y_n^0$, compute $r^0=\xi_{n}-A_ny_n^0$
\STATE Choose $\hat{r}=r^0$
\FOR{$k=1,2,\ldots,$}
\STATE $\rho^{k-1}=\hat{r}^Tr^{k-1}$
\IF{k=1}
\STATE $p^k=r^{k-1}$
\ELSE
\STATE $\beta^{k-1}=\left(\frac{\rho^{k-1}}{\rho^{k-2}}\right)\left(\frac{\alpha^{k-1}}{\omega^{k-1}}\right)$
\STATE $p^k=r^{k-1}+\beta^{k-1}\left(p^{k-1}-\omega^{k-1}v^{k-1}\right)$
\ENDIF
\STATE $\hat{p}=S_np^k, \hat{p}=S_n^T\hat{p}$
\STATE $v^k=A_n\hat{p}$
\STATE $\alpha^k=\frac{\rho^{k-1}}{\hat{r}^Tv^k}$
\STATE $s=r^{k-1}-\alpha^kv^k$
\IF{$||s||\le\epsilon$}
\STATE $y_n^k=y_n^{k-1}+\alpha^k\hat{p}$
\STATE  Stop
\ENDIF
\STATE $\hat{s}=S_ns, \hat{s}=S_n^T\hat{s}$
\STATE $t=A_n\hat{s}$
\STATE $\omega^k=\frac{t^Ts}{t^Tt}$
\STATE $y_n^k=y_n^{k-1}+\alpha^k\hat{p}+\omega^k\hat{s}$
\STATE $r^k=s-\omega^k t$
\IF{$||r^k||\le\epsilon$}
\STATE  Stop
\ENDIF
\ENDFOR
\end{algorithmic}
\end{algorithm}
\begin{algorithm}[!h t b]\label{Alg:3}
\caption{Calculating $A_ny_n$ by FFT}
\label{BICGStab}
\begin{algorithmic}[1]
\STATE Given initial vectors $a_1,q_1$
\STATE Set:\quad $c=[q_1,0,a_1,zeros(1, 2^n-2r+1)]^T$,  \,$t=y_n+\sqrt{-1}\hat{y_n}$\\[5pt]
\STATE \quad Do: $z=\mathcal{IF}(\mathcal{F}(t)\circ \mathcal{F}(c))$\\[5pt]
\STATE Set:\quad $A_ny_n:=ap\times real(z(1:2^n-r+1))+aq\times imag(z(2^n-r+1:-1:1))$\\[10pt]
\end{algorithmic}
\begin{algorithmic}
\STATE {\em Note:} $\mathcal{F},\mathcal{IF}$ denote the FFT and inverse FFT, respectively,  $\hat{y_n}=y_n(end:-1:1)$,\\
\qquad $\circ$ denotes the Hadamard product of vector $a$ and $b$.
\end{algorithmic}
\end{algorithm}

Now using the provided algorithms and taking $p=1$,~$q=0$, and $a=1$, we solve (\ref{1.1}) with $f(x)=f_1$ and $f_2$ given in (\ref{f1}) and (\ref{f2}), respectively, i.e.,
$$
\left\{
\begin{array}{c}
-D {}_0D_x^{-\beta} Du = f_{i}, \\ \\ 
u(0)=u(1) = 0,
\end{array}
\right.
$$
with $i=1,2$. When letting the exact solution be $u = x^{2}-x^{3}$, the forcing function is
\begin{equation} \label{f1}
f_{1} = \frac{-2x^{\beta}}{\Gamma(\beta+1)} + \frac{6x^{\beta+1}}{\Gamma(\beta+2)};
\end{equation}
and when the exact solution being taken as $u = x^{\lambda}-x$, we have the forcing function
\begin{equation}\label{f2}
f_{2} = \frac{-\Gamma(\lambda+1)x^{\lambda+\beta-2}}{\Gamma(\lambda+\beta-1)}+\frac{x^{\beta-1}}{\Gamma(\beta)}.
\end{equation}

The Bi-CGSTAB and the wavelet preconditioned Bi-CGSTAB methods are respectively used to solve the above equations. The numerical results are listed in Tables 3 and 4 for $f(x)=f_1$ and $f_2$, respectively. In performing the numerical computations, the initial value of iteration is taken as zero, and the stopping criterion $\epsilon=10^{-7}$. In fact, for making the comparisons, the Gaussian elimination based on the Doolittle LU decomposition (GE) is also used to solve the corresponding equations.
\begin{table}[!h t b]\fontsize{8pt}{12pt}\selectfont\label{tab:3}
\begin{center}
\caption{Numerical performances of the Bi-CGSTAB method, the GE method, and the preconditioned Bi-CGSTAB method, respectively, with the forcing function $f(x)=f_1$, $\beta=0.5$, and $r=2$ }\vspace{10pt}
\begin{tabular*}{\linewidth}{@{\extracolsep{\fill}}*{1}{c c c| c c c c}}
\hline
 $n$ & \multicolumn{2}{c|}{Bi-CGSTAB} & GE    &\multicolumn{2}{c}{Pre-Bi-CGSTAB}&$L_{2}$~error\\ \cline{5-6}
     &  iter     & cpu(s)       &cpu(s)&       iter  & cpu(s)          & \\
        \cline{2-3}                                  \cline{3-7}
5   &   33.5       & 0.0152        & 0.0073         &15.5   & 0.0132 &       2.3973e-04\\
6   & 70.5          &  0.0328       & 0.0277        &18.5   & 0.0180 &       6.0006e-05\\
7   &   142.5        & 0.0792       & 0.1264        &20.5   & 0.0241 &       1.5021e-05\\
8   &   312.5       & 0.1664        & 0.5030        &23.5   & 0.0334 &       3.7596e-06\\
9   &   783.5      & 0.5962       & 2.5780          &26.5   & 0.0583 &       9.4116e-07\\
10  &   1933.5       &  2.0709    & 14.9119        &27.0   & 0.1015 &       2.3875e-07\\ \hline
\end{tabular*}
\end{center}
\end{table}
\begin{table}[!h t b]\fontsize{8pt}{12pt}\selectfont
\begin{center}
\caption{Numerical performances of the Bi-CGSTAB method, the  GE method, and the preconditioned Bi-CGSTAB, respectively, with the forcing term $f(x)=f_2$, $\beta=0.75$, $\lambda=1.1$, and $r=2$}\vspace{10pt}\vspace{10pt}
\begin{tabular*}{\linewidth}{@{\extracolsep{\fill}}*{1}{c c c| c c c c}}
\hline
 $n$ &  \multicolumn{2}{c|}{Bi-CGSTAB}  &            GE    &\multicolumn{2}{c}{Pre-Bi-CGSTAB}&  $L_{2}$~error   \\ \cline{5-6}
     & iter   &  cpu(s)            &cpu(s)&       iter  & cpu(s)               & \\
        \cline{2-3} \cline{3-7}
5   &  54.5       & 0.0335            & 0.0135        &31   & 0.0230 &           1.0539e-04\\
6   &  134.5     & 0.0558             & 0.0466        &47   & 0.0370 &            3.4800e-05\\
7   &  246.5     & 0.1256             & 0.2065        &49   & 0.0466 &           1.1484e-05\\
8   &  516.5     & 0.2672            & 0.2845        &61   & 0.0648 &            3.7889e-06\\
9   &  1095.5     & 0.8060           & 4.3384        &65   & 0.1014 &            1.2499e-06\\
10  &  2471.5     & 2.5550         & 24.0359       &69   & 0.1736 &            4.1252e-07\\ \hline
\end{tabular*}
\end{center}
\end{table}
From Tables 3 and 4, it can be noted that both the Bi-CGSTAB method and the preconditioned Bi-CGSTAB method have a stable convergence rate $2$ or $1.6$ (due to the limited smoothness of the exact solution), but the preconditioned Bi-CGSTAB method is much faster than the Bi-CGSTAB method; the iteration numbers of the Bi-CGSTAB method increases quickly, but the ones of  the preconditioned Bi-CGSTAB method tends to be uniformly bounded.  It can also be noted that compared with the GE method the computational time can be greatly reduced, while they have almost the same $L^2$ error.

 For the sake of completeness, we also show the condition numbers of $A_{n}$ and $B_{n}$ with the spline basis of order $3$ in Tables 5 and 6.
\begin{table}[!h t b]\fontsize{8pt}{12pt}\selectfont
  \begin{center}
  \caption {The condition numbers of the matrixes $A_{n}$ and $B_{n}$ with $p = q = 0.5$, and $r = 3$} \vspace{5pt}
  \begin{tabular*}{\linewidth}{@{\extracolsep{\fill}}*{6}{c}}                                    \hline  
$\beta$ &$n$ & size of $B_{n}$  & $\kappa(A_{n})$        & $\kappa(B_{n})$ &      \\\hline
              &         3&    6 $\times$ 6&     4.6180&     6.0093  \\
              &         4&    14 $\times$ 14&     14.0369&      7.7079  \\
              &         5&    30 $\times$ 30&     41.3358&      9.2886  \\
$\beta = 0.5$ &         6&     62 $\times$ 62&     119.1580&     10.6787   \\
              &         7&     126 $\times$ 126&     340.1521&     11.8809  \\
              &         8&     254 $\times$ 254&     966.4690&      12.9127 \\
              &         9&     511 $\times$ 511&     2.7397e+003&     13.7796  \\\hline

               &        3&    6 $\times$ 6&     3.4303&     6.2506  \\
               &        4&    14 $\times$ 14&     8.8072&     7.6850  \\
               &        5&    30 $\times$ 30&     21.7487&      9.0615  \\
$\beta = 0.75$ &        6&     62 $\times$ 62&     52.6590&      10.2663   \\
               &        7&     126 $\times$ 126&     126.3300&     11.3041   \\
               &        8&     254 $\times$ 254&     301.7399&     12.1916 \\
               &        9&     510 $\times$ 510&     719.1686&      12.9483  \\\hline
    \end{tabular*}
  \end{center}
\end{table}
\begin{table}[!h t b]\fontsize{8pt}{12pt}\selectfont
  \begin{center}
  \caption {The condition numbers of the matrixes $A_{n}$ and $B_{n}$ with $p = 1$, $q = 0$, and $r = 3$} \vspace{5pt}
  \begin{tabular*}{\linewidth}{@{\extracolsep{\fill}}*{6}{c}}                                    \hline  
$\beta$ &$n$ & size of $B_{n}$  & $\kappa(A_{n})$        & $\kappa(B_{n})$ &      \\\hline
              &         3&    6 $\times$ 6&     4.9508&     6.0966  \\
              &         4&    14 $\times$ 14&     15.0241&      8.0872  \\
              &         5&    30 $\times$ 30&     43.9033&     9.8905  \\
$\beta = 0.5$ &         6&     62 $\times$ 62&     126.2291&      11.5143   \\
              &         7&     126 $\times$ 126&     359.9874&     12.9501   \\
              &         8&     254 $\times$ 254&     1.0224e+003&    14.2106 \\
              &         9&     511 $\times$ 511&     2.8980e+003&     15.2991  \\\hline

               &        3&    6 $\times$ 6&     4.4583&     6.9203  \\
               &        4&    14 $\times$ 14&     11.3404&      9.8872  \\
               &        5&    30 $\times$ 30&     26.7126&      12.0692  \\
$\beta = 0.75$ &        6&     62 $\times$ 62&     63.2405&      13.9770  \\
               &        7&     126 $\times$ 126&     150.4342&     15.7375   \\
               &        8&     254 $\times$ 254&     358.1752&     17.3488 \\
               &        9&     510 $\times$ 510&     852.6666&     18.8178  \\\hline
    \end{tabular*}
  \end{center}
\end{table}

\subsection{Wavelet Galerkin method for two dimensional FEDE}
 We know that the wavelet bases constructed in Lemma \ref{lemma:3.2} are the Riesz bases of $H_{0}^{\mu}((0,1)^2)$ with $0 < \mu < r - 1/2$; and they are applied to solve the following FEDE:
\begin{displaymath}
-D_{x}^{s}a_{1}(p_{1}{} _0D_x^{-\alpha} + q_{1}{} _x D_1^{-\alpha})D_{x}u - D_{y}^{s}a_{2}(p_{2}{} _0D_y^{-\beta} + q_{2}{} _y D_1^{-\beta})D_{y}u = f,~~x,y \in \Omega=[0,1]\times[0,1],
\end{displaymath}
where $D_{x}=\frac{\partial u(x,y)}{\partial x}$ and $D_{y}=\frac{\partial u(x,y)}{\partial y}$, and $D_{x~{\rm or}~y}^s$ denotes $s$-th derivative; ${} _0D_x^{-\alpha}$ (or ${} _0D_x^{-\beta}$) and ~${} _xD_1^{-\alpha}$ (or ${} _xD_1^{-\beta}$) represent the left and right Riemann-Liouville fractional integral operators, respectively, with $0 <\alpha,\beta <1$ satisfying $p_{1} + q_{1} = 1$ and $p_{2} + q_{2} = 1$. When $s = 2$ or $3$, the Riesz bases constructed in Lemma \ref{lemma:3.2} can be applied to solve this equation.

The variational formulation of this fractional differential equation is given as follows
\begin{eqnarray*}
B(u,v) &=& a_{1}p_{1}\langle D_{x}^{s-1}{} _0D_x^{-\alpha}D_{x}u,D_{x}v \rangle + a_{1}q_{1}\langle D_{x}^{s-1}{} _x D_1^{-\alpha}D_{x}u,D_{x}v \rangle \\
&& + a_{2}p_{2}\langle D_{y}^{s-1}{} _0D_y^{-\beta}D_{y}u,D_{y}v \rangle + a_{2}q_{2}\langle D_{y}^{s-1}{} _y D_1^{-\beta}D_{y}u,D_{y}v \rangle.
\end{eqnarray*}

Consequently, in order to solve the variational formulation in $\Omega=(0,1)^2$, we use  $\widetilde{V}_{k}$  to approximate the $H_{0}^{\mu}(\Omega)$ space, since $\widetilde{\Phi}_{n} := \{\widetilde{\phi}_{n,(j_{1},j_{2})}:j_{1}=0,1,\ldots,2^{n}-r;\,j_{2}=0,1,\ldots,2^{n}-r\}$ is a basis of $\widetilde{V}_{n}$. We investigate a $u_{k} \in \widetilde{V}_{n}$ such that
\begin{eqnarray}\label{equ:6.1.3}
&&{} a_{1}p_{1}\langle D_{x}^{s-1}{} _0D_x^{-\alpha}D_{x}u_{k},D_{x}v \rangle + a_{1}q_{1}\langle D_{x}^{s-1}{} _x D_1^{-\alpha}D_{x}u_{k},D_{x}v \rangle +\\
&&{} a_{2}p_{2}\langle D_{y}^{s-1}{} _0D_y^{-\beta}D_{y}u_{k},D_{y}v \rangle + a_{2}q_{2}\langle D_{y}^{s-1}{} _y D_1^{-\beta}D_{y}u_{k},D_{y}v \rangle = \langle f,v \rangle ,  ~~\forall v \in \widetilde{V}_{n}.\nonumber
\end{eqnarray}
Suppose $u_{n} = \sum\nolimits_{\tilde{\phi} \in \widetilde{\Phi}_{n}}y_{\tilde{\phi}}\tilde{\phi}$. Let~$C_{k}^{s}$~be the matrix $(a_{1}p_{1}\langle D_{x}^{s-1}{} _0D_x^{-\alpha}D_{x}\tilde{\sigma},D_{x}\tilde{\phi} \rangle + a_{1}q_{1}\langle D_{x}^{s-1}{} _xD_1^{-\alpha}D_{x}\tilde{\sigma},D_{x}\tilde{\phi} \rangle \\+ a_{2}p_{2}\langle D_{y}^{s-1}{} _0D_y^{-\beta}D_{y}\tilde{\sigma},D_{y}\tilde{\phi} \rangle +  a_{2}q_{2}\langle D_{y}^{s-1}{} _yD_1^{-\beta}D_{y}\tilde{\sigma},D_{y}\tilde{\phi} \rangle)_{\tilde{\sigma},\tilde{\phi} \in \widetilde{\Phi}_{n}}$; and $\widetilde{\xi}_{n}$ be the column vector $(\langle f,\tilde{\phi} \rangle)_{\tilde{\phi} \in \widetilde{\Phi}_{n}}$.  Then the column vector $y_{n}=(y_{\tilde{\phi}})_{\tilde{\phi} \in \widetilde{\Phi}_{n}}$ is the solution of the linear system
\begin{equation}\label{equ:6.1.4}
C_{n}^{s}y_{n} = \widetilde{\xi}_{n}.
\end{equation}
Similar to the one dimensional case, without preconditioning it would be difficult to solve the system when we increase the discrete level  of $n$.

Now we employ the wavelet bases constructed above to solve the variational problem. For $n \geq n_{0}$, $\widetilde{\Psi}_{n} := \{2^{-n_{0}\mu}\widetilde{\phi}_{n_{0},j}:j \in J_{n_{0}}\} \cup \bigcup_{k=k_{0}}^{n-1}\{2^{-k\mu}\omega:\omega \in \Gamma_{k}\}$. To find a column vector $z_{n}=(z_{\widetilde{\psi}})_{\widetilde{\psi} \in \widetilde{\Psi}_{n}}$ such that
\begin{eqnarray}\label{equ:6.1.5}
D_{n}^{s}z_{n} = \widetilde{\eta}_{n},
\end{eqnarray}
where $D_{n}^{s}$ is matrix $(a_{1}p_{1}\langle D_{x}^{s-1}{} _0D_x^{-\alpha}D_{x}\widetilde{\chi},D_{x}\widetilde{\psi} \rangle +a_{1}q_{1}\langle D_{x}^{s-1}{} _xD_1^{-\alpha}D_{x}\widetilde{\chi},D_{x}\widetilde{\psi} \rangle)_{\widetilde{\chi},\widetilde{\psi} \in \widetilde{\Psi}_{n}}+(a_{2}p_{2}\langle D_{y}^{s-1}{} _0D_y^{-\beta}\\D_{y}\widetilde{\chi},D_{y}\widetilde{\psi} \rangle + a_{2}q_{2}\langle D_{y}^{s-1}{} _yD_1^{-\beta}D_{y}\widetilde{\chi},D_{y}\widetilde{\psi} \rangle)_{\widetilde{\chi},\widetilde{\psi} \in \widetilde{\Psi}_{n}}$, and $\widetilde{\eta}_{k}$ denotes the column vector $(\langle f,\widetilde{\psi} \rangle)_{\widetilde{\psi} \in \widetilde{\Psi}_{}}$. Then  the approximate solution in $\widetilde{V}_{n}$ can be written as $u_{n}=\sum_{\widetilde{\psi} \in \widetilde{\Psi}_{n}}z_{\widetilde{\psi}}\widetilde{\psi}$.

Since $\widetilde{\Phi}_{n}$ and $\widetilde{\Psi}_{n}$ are two different bases of $\widetilde{V}_{n}$, there is a unique transformation~$\widetilde{S}_{n}$ between the two bases. Consequently, we have
\begin{displaymath}
D_{n}^{s} = \widetilde{S}_{n}C_{n}^{s}\widetilde{S}_{n}^{T},~~\widetilde{\eta}_{n} = \widetilde{S}_{n}\widetilde{\xi}_{n},
\end{displaymath}
and \eqref{equ:6.1.4} is equivalent to \eqref{equ:6.1.5}. If we set $y_{n} = \widetilde{S}_{n}^{T}z_{n}$, then \eqref{equ:6.1.3} is preconditioned by the matrix $\tilde{S}_{n}$.

Next, we provide the algorithm of generating $D_{n}^{s}$, where $s = 2$ or $s = 3$. For generating the matrix $D_n^s$, we need the transform matrix $\widetilde{S}_{n}$, which is far more complex than one dimensional case; the one-level transform is given in Algorithm 4, and the others are like the steps proposed in Algorithm 1, which are omitted here.

\begin{algorithm}[!h t b]
\caption{Generating the matrix $\tilde{P_k}$ for $2D$}
\begin{algorithmic}[1]

\STATE Let $P_{k}$ be one-level transform matrix in $1D$, satisfying:
      ${\Phi_{k-1} \choose 2^{-(k-1)\mu}\Gamma_{k-1}} = P_k\Phi_{k}$
\STATE Initialize the matrixes:\\
      $t_1=2^{k-1}-r+1$, $t_2=2^{k}-r+1$

      $L1=a_1\otimes P_k(1:t_1,1:t_2)$\\
      $L2=a_1\otimes P_k(t_1+1:t_2,1:t_2)$\\

      $L3=a_2\otimes P_k(1:t_1,1:t_2)$\\
      $L4=a_3\otimes P_k(1:t_1,1:t_2)$\\
      $L5=\hat{a}_2 \otimes P_k(1:t_1,1:t_2)$\\

      $L6=a_2\otimes P_k(t_1+1:t_2,1:t_2)$\\
      $L7=a_3\otimes P_k(t_1+1:t_2,1:t_2)$\\
      $L8=\hat{a}_2\otimes P_k(t_1+1:t_2,1:t_2)$
\FOR{$i=1:2^{k-1}-r+1$}
  \STATE $A=\left[T_bA;T_f^{(i-1)} L_1\right]$, \quad $C=\left[T_bC;T_f^{(i-1)}L_2\right]$
\ENDFOR

\STATE $L_3=\left[L_3,zeros\left(2^{k-1}-r+1,(2^k-r+1)^2-len_{col}(L_3)\right)\right]$
\STATE $L_6=\left[L_6,zeros\left(2^{k-1},(2^k-r+1)^2-len_{col}(L_6)\right)\right]$
\FOR{$i=1:2^{k-1}-2$}
\STATE $B=\left[T_bB,T_f^{(i-1)}L_4\right]$,\quad $D=\left[T_bD,T_f^{(i-1)}L_7\right]$
\ENDFOR
\STATE $L_5=\left[zeros\left(2^{k-1}-r+1,(2^k-r+1)^2-len_{col}(L_5)\right),L_5\right]$,
\STATE $L_8=\left[zeros\left(2^{k-1},(2^k-r+1)^2-len_{col}(L_8)\right),L_8\right]$

\STATE $\tilde{P_k}\left(1:(2^{k-1}-r+1)^{2}, 1:(2^{k}-r+1)^{2}\right) = A$
\STATE $\tilde{P_k}\left((2^{k-1}-r+1)(2^k-r+1)+1:(2^{k-1}-r+1)(3\times 2^{k-1}-r+1), 1:(2^{k}-r+1)^{2}\right) = C$
\STATE $\tilde{P_k}\left(2^{k-1}-r+1)^2+1:(2^{k-1}-r+1)(2^k-r+1), 1:(2^{k}-r+1)^{2}\right) = \left[L_3;\,B;\,L_5\right]$
\STATE $\tilde{P_k}\left((2^{k-1}-r+1)(3\times 2^{k-1}-r+1)+1:(2^{k}-r+1)^2, 1:(2^{k}-r+1)^{2}\right) = \left[L_6;\,D;\,L_8\right]$\\[9pt]
\end{algorithmic}
\begin{algorithmic}
\STATE {\em Note:} $a_1,a_2,a_3$ denote the refinement coefficient vectors.\\
For example, when $r=3$, we have:\\
  $a_1=[\frac{1}{4},\frac{3}{4},\frac{3}{4},\frac{1}{4}],\quad a_2=[\frac{5}{12},-\frac{1}{12}]$,\\
  $a_3=[\frac{1}{12},-\frac{5}{12},\frac{5}{12},-\frac{1}{12}],\quad\hat{a}_2=a_2(end:-1:1)$.\\
  $a\otimes b$ denotes the kronecker product of vector $a$ and $b$; and $len_{col}(L)$ denotes the column number of matrix $L$.\\
  $For \left[T_bG;T_f^{(i-1)} L\right]$, $T_b$ denotes the zero padding operator for making the column of $T_bG$ equal to $T_f^{(i-1)}L$; for simplicity, $A,B,C,D $ are initialized with empty matrix.
  $T_f$ is the operator to extend the matrix $L$ by adding new columns with the value of zero at the left hand side of the matrix. More precisely, for $L=L_1$ and $L4$, $T_fL$ adds $2^k-2r+2$ column zeros before $L$, but for $L=L_2$ and $L_7$, it adds $2^k$ column zeros before $L$.
\end{algorithmic}
\end{algorithm}

Taking $u=x^{2}(1-x)^{2}y^{2}(1-y)^{2}, \alpha=\beta=\frac{3}{4}, p_1=p_2=1, q_1=q_2=0, a_1=a_2=0$,  then $u$ is the exact solution of the following equation
\begin{equation*}
-D_{x}^{s}{}_0D_x^{-3/4}Du - D_{y}^{s}{}_0D_y^{-3/4}Du = f_{s}\\[4pt]
\end{equation*}
with the boundary conditions  $u(x,y)|_{\partial \Omega}=0$, $(\partial u(x,y)/ \partial x)|_{x=0,y\in [0,1]}=0$, and $(\partial u(x,y)/ \partial y)|_{x\in [0,1],y=0}=0$ for $s=2$,  and $u(x,y)|_{\partial \Omega}=0$, $(\partial u(x,y)/ \partial x)|_{x=0~{\rm and}~ x=1,\,y\in [0,1]}=0$, and $(\partial u(x,y)/ \partial y)|_{x\in [0,1],\,y=0~ {\rm and} ~y=1}=0$ for $s=3$, respectively; and
the forcing function $f_s$  is given as follows
\begin{eqnarray*}
f_s&=&\left(\frac{-24}{\Gamma(\frac{19}{4}-s)}x^{\frac{15}{4}-s}+\frac{12}{\Gamma(\frac{15}{4}-s)}x^{\frac{11}{4}-s}+
\frac{-2}{\Gamma{(\frac{11}{4}-s)}}x^{\frac{7}{4}-s}\right)y^2(1-y^2)\\
  &&+\left(\frac{-24}{\Gamma(\frac{19}{4}-s)}y^{\frac{15}{4}-s}+\frac{12}{\Gamma(\frac{15}{4}-s)}y^{\frac{11}{4}-s}+
\frac{-2}{\Gamma{(\frac{11}{4}-s)}}y^{\frac{7}{4}-s}\right)x^2(1-x^2).
\end{eqnarray*}

We first calculate the condition numbers of the corresponding stiffness matrixes, then use the algorithms presented in previous sections to compute the numerical solutions for different $s$ with $r=3$. For confirming the relation between the conditional numbers and $\beta$, we also list the condition numbers of the matrixes with $\beta=0.25$.
\begin{table}[!h t b]\fontsize{8pt}{12pt}\selectfont
  \begin{center}
  \caption {When $r = 3$, $\alpha=\beta=0.75$, the condition number of the matrix $C_{n}^{s}$ and $D_{n}^{s}$} \vspace{5pt}
  \begin{tabular*}{\linewidth}{@{\extracolsep{\fill}}*{7}{c}}                                    \hline
   $s$      &  $n$  &  size of $D_{n}^{s}$   & $\kappa(C_{n}^{s})$   &ratio       & $\kappa(D_{n}^{s})$  & ratio    \\\hline
            &   4   &    196$\times$ 196     &    29.4463            &               &54.8262              &           \\
   $s=2$    &   5   &    900$\times$ 900     &   144.4390            &2.2943         &59.8263              & 0.1206     \\
            &   6   &   3844$\times$3844     &   698.5068            &2.2738         &63.6789              & 0.0954      \\\hline

            &   4   &    196$\times$ 196     &   160.09e+02           &               &82.9191               &           \\
   $s=3$    &   5   &    900$\times$ 900     &   1.5444e+03           &3.2701         &107.6680              & 0.3768     \\
            &   6   &    3844$\times$3844    &   1.4752e+04           &3.2558         &116.0167              & 0.1077      \\\hline
    \end{tabular*}
  \end{center}
\end{table}
\begin{table}[!h t b]\fontsize{8pt}{12pt}\selectfont
  \begin{center}
  \caption {When $r = 3$, $\alpha=\beta=0.25$, the condition number of the matrix $C_{n}^{s}$ and $D_{n}^{s}$} \vspace{5pt}
  \begin{tabular*}{\linewidth}{@{\extracolsep{\fill}}*{7}{c}}                                    \hline
   $s$      &  $n$  &  size of $D_{n}^{s}$   & $\kappa(C_{n}^{s})$   &ratio       & $\kappa(D_{n}^{s})$  & ratio    \\\hline
            &   4   &    196$\times$ 196     &   7.4192e+01          &            &69.8276               &           \\
   $s=2$    &   5   &    900$\times$ 900     &   5.1254e+02          &2.7883      &89.9318               & 0.3650     \\
            &   6   &   3844$\times$3844     &   3.4862e+03          &2.7659      &94.9550               & 0.0784     \\\hline

            &   4   &    196$\times$ 196     &   5.8919e+02           &              &159.0936             &           \\
   $s=3$    &   5   &    900$\times$ 900     &   8.0236e+03           &3.7675        &176.1326             & 0.1468    \\
            &   6   &    3844$\times$3844    &   1.0824e+05           &3.7539        &184.4045             & 0.0662     \\\hline
    \end{tabular*}
  \end{center}
\end{table}


From Tables 7 and 8, it can be noted that the condition numbers of the stiffness matrix corresponding to the ordinary Galerkin methods increase with a rate as  $\mathcal{O}(h^{-(s+1-\beta)})$, but the conditional numbers corresponding to the wavelet Galerkin methods tend to be uniformly bounded, as stated in Lemma \ref{lemma:3.2}. For the numerical iterative schemes, because of the tensor form of the matrixes $C_{n}^{s}$, we can still make use of the Toeplitz structure of the matrix and FWT such that the computation complexity is as $\mathcal{O}(N\log N)$, where $N$ denotes the number of bases. The numerical performances for $\beta=0.75$ are presented in Table 9, 10, 11, and 12. It can be seen that for getting the same accuracy, compared with the Bi-CGSTAB method, the preconditioned Bi-CGSTAB method needs less computation time and the number of iterations when $s=2$ and $3$; and in fact when $s=3$, the numerical errors for Bi-CGSTAB method increase early. For the GMRES method and the preconditioned GMRES method, they have almost the same $L^2$ errors, but the latter method converges much faster. From Tables 11 and 12, it can be noted that GMRES(50) is faster than GMRES(20), but preconditioned GMRES(20) is faster than preconditioned GMRES(50).

\begin{table}[!h t b]\fontsize{8pt}{12pt}\selectfont
\begin{center}
\caption{Numerical performances for the Bi-CGSTAB method and the preconditioned Bi-CGSTAB method with $\alpha=\beta=0.75$, and $ r=3$} \vspace{10pt}
\begin{tabular*}{\linewidth}{@{\extracolsep{\fill}}*{1}{c  c| c c c |c c c}}
\hline
 $s$ &  $n$ & \multicolumn{3}{c|}{Bi-CGSTAB} &\multicolumn{3}{c}{Pre-Bi-CGSTAB}   \\
 &          & iter           &cpu(s)            &$L_{2}$~error    &iter    & cpu(s) &$L_{2}$~error             \\\hline
      &5    & 43.0            & 0.2046          & 6.1592e-07     &41.5     &0.2315 &     6.1592e-07\\
 $s=2$& 6   & 109.5          & 1.3395         & 8.8556e-08     &46.5     & 0.8530 &     8.8556e-08\\
      &7    & 277.5          & 8.3340         &1.3326e-08      &55.5     & 3.0553 &     1.3326e-08 \\\hline

     &5    & 347.5           & 1.5635         &5.8094e-06      &348      & 1.8908 &      5.8094e-06\\
 $s=3$  &6 & 996             & 15.0287        &8.1043e-05      &403      & 7.8870 &      1.4475e-06 \\
     &7    & 2569            & 88.5789        &1.3419e-03      &469      & 25.5683 &     3.6123e-07 \\ \hline

\end{tabular*}
\end{center}
\end{table}
\begin{table}[!h t b]\fontsize{8pt}{12pt}\selectfont
\begin{center}
\caption{Numerical performances for the non-restarted GMRES method and the preconditioned GMRES method with $\alpha=\beta=0.75$, and $ r=3$} \vspace{10pt}
\begin{tabular*}{\linewidth}{@{\extracolsep{\fill}}*{1}{c  c| c c |c c| c}}
\hline
 $s$ &  $n$ & \multicolumn{2}{c|}{GMRES} &\multicolumn{2}{c|}{Pre-GMRES}  &$L_{2}$~error \\
 &          & iter           &cpu(s)             &iter     & cpu(s)   &              \\\hline
      &5    & 60           & 0.3359             &55        &1.3109 &            6.1592e-07\\
 $s=2$& 6   & 135          & 3.3185             &61        & 1.1250&             8.8556e-08\\
      &7    & 309          & 47.443             &65        & 3.9323 &            1.3326e-08 \\\hline

     &5    & 231           & 3.0238             &185       & 2.0874 &            5.8094e-06\\
 $s=3$  &6 & 721           & 70.1586            &260       & 7.8870 &            1.4475e-06 \\
     &7    & 2670          & 3299.2131          &313       & 53.7952&            3.6123e-07 \\ \hline

\end{tabular*}
\end{center}
\end{table}
\begin{table}[!h t b]\fontsize{8pt}{12pt}\selectfont
\begin{center}
\caption{Numerical performances for the restarted GMRES(50) method and the preconditioned GMRES(50) method with $\alpha=\beta=0.75$, and $ r=3$} \vspace{10pt}
\begin{tabular*}{\linewidth}{@{\extracolsep{\fill}}*{1}{c  c| c c |c c |c}}
\hline
 $s$ &  $n$ & \multicolumn{2}{c|}{GMRES(50)} &\multicolumn{2}{c|}{Pre-GMRES(50)}  &$L_{2}$~error \\
 &          & iter               &cpu(s)             &iter     & cpu(s)   &              \\\hline
     &5    & 1$\times$50+15           & 0.2746            &1$\times$50+7       & 0.2838 &           6.1592e-07 \\
 $s=2$  &6 & 6$\times$50+5            & 3.7541             &1$\times$50+14      & 0.9887 &          8.8556e-08 \\
     &7    & 14$\times$50+35          & 24.8405            &1$\times$50+16       & 3.3791&            1.3343e-08 \\ \hline
       &5    & 20$\times$50+5        & 4.4621e+00              &4$\times$50+27       &1.1295 &        5.8094e-06       \\
 $s=3$& 6   & 93$\times$50+32        & 5.5601e+01             &5$\times$50+28       & 4.2065&             1.4475e-06\\
      &7    & 1478$\times$50+12         &2.5654e+03             &6$\times$50+47       & 16.1958 &             3.6123e-07 \\\hline
\end{tabular*}
\end{center}
\end{table}
\begin{table}[!h t b]\fontsize{8pt}{12pt}\selectfont
\begin{center}
\caption{Numerical performance with the restarted GMRES(20) method and the preconditioned GMRES(20) method, respectively, $\alpha=\beta=0.75$, $ r=3$.} \vspace{10pt}
\begin{tabular*}{\linewidth}{@{\extracolsep{\fill}}*{1}{c  c| c c |c c |c}}
\hline
 $s$ &  $n$ & \multicolumn{2}{c|}{GMRES(20)} &\multicolumn{2}{c|}{Pre-GMRES(20)}  &$L_{2}$~error \\
 &          & iter               &cpu(s)             &iter     & cpu(s)   &              \\\hline
     &5    & 6$\times$20+14           & 0.4269            &13$\times$20+14       &  0.2193 &           6.1592e-07 \\
 $s=2$  &6 & 17$\times$20+8            & 3.0920             &17$\times$20+6       & 0.7959 &          8.8556e-08 \\
     &7    & 81$\times$20+12          & 40.0962            &20$\times$20+14       & 2.5415&            1.3343e-08 \\ \hline
       &5    & 65$\times$20+16        & 4.4621e+00              &13$\times$20+14      &1.0174 &        5.8094e-06       \\
 $s=3$& 6   & 455$\times$20+16        & 5.5601e+01             &18$\times$20+6       &  4.2484&             1.4475e-06\\
      &7    & 12125$\times$20+10         &2.5654e+03             &20$\times$20+14       & 14.6722 &             3.6123e-07 \\\hline
\end{tabular*}
\end{center}
\end{table}

\section{Conclusion}
For improving the efficiency of solving FEDEs, three natural ways can be adopted: 1. reducing matrix vector multiplication from $\mathcal{O}(N^2)$ to $\mathcal{O}(N\log N)$; 2. keeping the condition numbers small and uniformly bounded; 3. increasing the convergence orders. For the general linear finite element methods,  because of the potential Toeplitz structure of the stiffness matrix, the cost of the matrix vector multiplication can be kept as $\mathcal{O}(N\log N)$. But for the high order elements, the potential Toeplitz structure is destroyed and the cost of matrix vector multiplication is $\mathcal{O}(N^2)$. If taking the scale functions (generally used to generate wavelets) as the base functions of the Galerkin methods, the potential Toeplitz structure of the stiffness matrix can be kept. Furthermore, based on the general wavelet theory, the Riesz bases of the space that the FEDE works are found and effectively used to solve the one and two dimensional FEDEs. The detailed algorithm descriptions are presented. The extensive numerical experiments are performed, and the numerical observables, including the condition numbers, iteration numbers, cpu time cost, are calculated; all demonstrate the striking benefits of the wavelet Galerkin methods in solving FEDEs.


\section*{Acknowledgements}
This work was supported by the National Natural Science Foundation of China under Grant  No. 11271173.

\end{document}